\documentclass[letterpaper, 10 pt, conference, onecolumn]{ieeeconf}  

\IEEEoverridecommandlockouts                              

\usepackage[utf8x]{inputenc}

\usepackage[greek,english]{babel}
\usepackage{gfsporson}

\usepackage{balance}
\usepackage{graphicx}
\usepackage{siunitx}
\usepackage{cuted} 
\usepackage{amsmath} %
\usepackage{amsfonts}
\usepackage{amssymb}  %
\usepackage{mathrsfs}
\usepackage{mathtools}
\usepackage{url}

\makeatletter
\let\NAT@parse\undefined
\makeatother

\let\ieeebibliography\thebibliography

\usepackage[numbers,sort&compress]{natbib}

\renewcommand\thebibliography[1]{\ieeebibliography{#1}}

\usepackage[colorlinks]{hyperref}
\usepackage[authormarkup=none,commentmarkup=uwave,final]{changes}
\definechangesauthor[color=purple]{TC}
\definechangesauthor[color=violet]{IK}

\usepackage[IEEE]{altthm-styles}

\DeclareMathOperator{\distance}{dist}
\DeclareMathOperator{\kernel}{ker}

\title{\LARGE \bf
\added[id=TC]{Vector Lyapunov functions for nonconvex constrained optimization} 
}

\author{Torbjørn Cunis$^{1}$ and Ilya Kolmanovsky$^{2}$
\thanks{$^{1}$TC is with the Institute of Flight Mechanics and Controls, University of Stuttgart, 70569 Stuttgart, Germany, and with the Department of Aerospace Engineering, University of Michigan, Ann Arbor, MI 48109, USA.
        \url{tcunis@ifr.uni-stuttgart.de}, \url{tcunis@umich.edu}}%
\thanks{$^{2}$IK is with the Department of Aerospace Engineering, University of Michigan, Ann Arbor, MI 48109, USA.
        \url{ilya@umich.edu}}%
}

\begin{document}

\maketitle
\thispagestyle{empty}
\pagestyle{empty}

\begin{abstract}
Nonconvex constrained optimization problems are commonly solved by \added[id=TC]{a Newton method for necessary conditions (sequential linear/quadratic programming) or a barrier reformulation (interior-point methods)}.
We present a vector Lyapunov analysis for Newton-type algorithms based on constraint violation and optimality of the primal solutions to a sequence of approximate optimization problems. We thus prove local stability and convergence by constructing a strict Lyapunov function for the primal solution.
\end{abstract}

\section{Introduction}
Lyapunov function arguments are strongly connected to the objectives of optimization algorithms \cite{polyak_lyapunov_2017}. Recent works have been dedicated to the construction of Lyapunov functions for (accelerated) gradient descent methods \cite{laborde_lyapunov_2020, wilson_lyapunov_2021, sanz_serna_connections_2021}. Conversely, globalization techniques for unconstrained optimization maintain `sufficient' decrease in the cost function \cite{Nocedal2006}.

In constrained optimization, however, the search for a feasible solution conflicts with the objective of minimization. Penalty methods combine cost function and constraint violation into {\em merit} functions, but require careful adjustment of the penalty parameter. 
Works \cite{dhingra_proximal_2019, mai_convergence_2020, hassan-moghaddam_proximal_2021} constructed Lyapunov functions for projected or proximal gradient descent with convex constraints. Similarly, \cite{cherukuri_asymptotic_2016, cherukuri_role_2018, qu_exponential_2019} studied primal-dual dynamics for saddle-point problems in convex constrained optimization, and \cite{boyd_distributed_2010} constructed a Lyapunov function for a primal-dual solution of the Alternating Direction of Multiplier Method with linear constraints. Other authors applied feedback stabilization to nonlinear programs \cite{karafyllis_feedback_2014, karafyllis_global_2017, allibhoy_control-barrier-function-based_2024}: While \cite{allibhoy_control-barrier-function-based_2024} employed control barrier functions to ensure that constraints are satisfied, \cite{karafyllis_global_2017} proposed two Lyapunov-like functions -- one for constraint violation, one for the cost -- and designed a dynamical system which reduces constraint violation and minimizes the cost on the feasible set. This last result is reminiscent of the vector Lyapunov method proposed by Matrosov \cite{matrosov_vector_1997}.

This work is concerned with a Lyapunov-based analysis of the {\em primal} solution in Newton-type algorithms for constrained nonconvex optimization problems, such as sequential convex (linear or quadratic) programming or interior-point methods. Such algorithms have eluded the construction of Lyapunov functions, including from Matrosov conditions, since Newton's method {\em jointly} solves for optimality and feasibility. Instead, Newton-type methods for optimization are commonly analyzed in terms of the primal-dual solutions to the necessary conditions \cite{Coleman1984, Dokov1998, Artacho2011, Colabufo2020}. Although Newton's method (for necessary conditions) suggests a Lyapunov candidate that combine constraint violation and first-order stationarity, such an approach neglects the objective of solving the {\em optimization} problem.

Instead, we are going to propose a pair of Lyapunov-like value functions for constraint violation and relative optimality, respectively. We then prove that a Newton iteration generated by a sequence of approximate quadratic optimization problems with linearized constraints (locally) satisfies a set of dissipation inequalities similar to those in network dynamics \cite{Arcak2016} and input-to-state stability \cite{Sontag2008}. Finally, the construction of a strict Lyapunov function, inspired by \cite{malisoff_constructions_2009}, certifies local stability of the primal solution with linear (or better) convergence. For the ease of notation, we restrict this work to equality constraints. Our analysis is based on three assumptions on the approximate problems: surjectivity of the linear constraints (ensures feasibility), strong convexity of the quadratic cost function (ensures unique solutions), and a `sufficiently' good approximation of the nonconvex necessary conditions. Similar assumptions, albeit formulated for the necessary conditions, have been employed in previous works on Newton-type methods (see, e.g., \cite{Cunis2024cdc}). Assumptions of convexity and surjectivity were made also in \cite{wachter_line_2005} for convergence of a filter-based globalization. \added[id=TC]{Our work facilitates a Lyapunov-based analysis of Newton-type algorithms for nonconvex optimization and their interaction with dynamical systems.}


\paragraph*{Notation}
For a function $f: \mathbb R^n \to \mathbb R^m$ and $x \in \mathbb R^n$, the {\em strict Fréchet derivative} of $f$ at $x$ is the linear mapping $\nabla f(x): \mathbb R^n \to \mathbb R^m$ satisfying
\begin{align*}
    \limsup_{\substack{x_1, x_2 \to x \\ x_1 \neq x_2}} \frac{\| f(x_1) - \nabla f(x) (x_1 - x_2) - f(x_2) \|}{\| x_1 - x_2 \|} = 0
\end{align*}
If such a mapping exist, $f$ is called {\em strictly Fréchet differentiable} at $x$ and $\nabla f(x)$ is unique. For a linear mapping $A: \mathbb R^n \to \mathbb R^m$, the {\em adjoint} of $A$ is the (linear) mapping $A^*: \mathbb R^m \to \mathbb R^n$ satisfying $\langle y, A x \rangle = \langle A^* y, x \rangle$ for all $x \in \mathbb R^n$ and $y \in \mathbb R^m$, where $\langle \cdot, \cdot \rangle$ denotes the inner product.

\section{Problem Statement}
We consider an equality-constrained nonlinear optimization problem of the form,
\begin{align}
    \label{eq:problem}
    \min_{x \in \mathbb R^n} \varphi(x) \quad \text{subject to $g(x) = 0$,}
\end{align}
where the cost function $\varphi: \mathbb R^n \to \mathbb R$ and constraint function $g: \mathbb R^n \to \mathbb R^m$ are strictly Fréchet differentiable \added[id=IK]{locally around a (local) optimal solution $\bar x \in \mathbb R^n$}. The feasible set of \eqref{eq:problem} is given by
\begin{align*}
    \Omega = \{ x \in \mathbb R^n \mid g(x) = 0 \}
\end{align*}
and the Lagrange function is 
\begin{align*}
    L: (x,y) \mapsto \varphi(x) + \langle y, g(x) \rangle
\end{align*}
with dual variable $y \in \mathbb R^m$.

For a given solution iterate, $x_k \in \mathbb R^n$, $k \in \mathbb N$, we approximate \eqref{eq:problem} by an optimization problem of minimizing a quadratic cost function subject to linear equality constraints:
\begin{align}
    \label{eq:quadratic-problem}
\begin{aligned}
    &\min_{x \in \mathbb R^n} \varphi_k(x) := \langle x - x_k, H_k (x - x_k) \rangle + \nabla \varphi(x_k)(x - x_k) \\
    &\text{subject to $g_k(x) := g(x_k) + \nabla g(x_k)(x - x_k) = 0$}
\end{aligned}
\end{align}
where $H_k: \mathbb R^n \to \mathbb R^n$ is a symmetric matrix. 
For (\ref{eq:quadratic-problem}), we define the feasible set as
\begin{align*}
    \Omega_k = \{ x \in \mathbb R^n \mid g_k(x) = 0 \}
\end{align*}
and the Lagrange function as
\begin{align*}
    L_k: (x,y) \mapsto \varphi_k(x) + \langle y, g_k(x) \rangle
\end{align*}
with $y \in \mathbb R^m$. 

If \eqref{eq:quadratic-problem} has an optimal solution, it becomes the next solution iterate. We denote the solution map of \eqref{eq:quadratic-problem}, parametrized in $x_k$, by $\Phi_k: \mathbb R^n \rightrightarrows \mathbb R^n$; and \added[id=TC]{its {\em dual} solutions by}
\begin{align*}
    \Psi_k(x_k) = \{ y \in \mathbb R^m \mid \nabla L_k(x, y) = 0, x \in \Phi_k(x_k) \}.
\end{align*}
We are interested in the convergence of a sequence $\{ x_k \}_{k \in \mathbb N}$ generated by the iteration
\begin{align}
    \label{eq:newton-iteration}
    x_{k+1} \in \Phi_k(x_k), \quad k \in \mathbb N,
\end{align}
existence provided, in the neighbourhood of \added[id=TC]{$\bar x$}. 

\begin{assumption}[Standing assumptions]
    \label{ass:standing}
    There exists an open neighbourhood $U \subset \mathbb R^n$ of $\bar x$ such that, for all $x_k \in U$ and $k \in \mathbb N$, the quadratic optimization \eqref{eq:quadratic-problem} satisfies the following properties:

    \paragraph{Regularity of $g_k$} The derivative $\nabla g_k(x_k): \mathbb R^n \to \mathbb R^m$ is surjective;
    \paragraph{Strong convexity of $\varphi_k$ relative to $\Omega_k$} There exists a constant $\mu > 0$ such that
    \begin{align*}
        \varphi_k(x_1) \geq \varphi_k(x_2) + \langle \nabla \varphi_k(x_2) +  \added{\nabla g_k(x_2)^* y}, x_1 - x_2 \rangle + \frac{\mu}{2} \| x_1 - x_2 \|^2
    \end{align*}
    for all $x_1, x_2 \in \Omega_k$ and \added{$y \in \mathbb R^m$};
    \paragraph{Lipschitz continuity of $\nabla[L - L_k](\cdot, y)$} There exists a constant $\gamma > 0$ such that\footnotemark
    \begin{align*}
        \| \nabla [L - L_k](x_1, \added[id=TC]{y_k}) - \nabla [L - L_k](x_2, \added[id=TC]{y_k}) \| \leq \gamma \| x_1 - x_2 \|
    \end{align*}
    for all $x_1, x_2 \in U$ and \added[id=TC]{$y_k \in \Psi_k(x_k)$};
    and $\mu > \gamma$. 
\end{assumption}

\footnotetext{Throughout this paper, $\nabla L(\cdot, y)$ and $\nabla L_k(\cdot, y)$ denote the strict Fréchet derivative of $L(\cdot, y)$ and $L_k(\cdot, y)$, respectively, for a fixed $y \in \mathbb R^m$.}

The regularity condition is a common constraint qualification; if $\nabla g_k(x_k)$ is surjective, an optimal solution $\bar x_k \in \mathbb R^n$ of \eqref{eq:quadratic-problem} necessarily satisfies
\begin{align}
    \label{eq:stationarity-quadratic}
    \nabla \varphi_k(\bar x_k) + \nabla g_k(\bar x_k)^* y_k = 0
\end{align}
for some $y_k \in \mathbb R^m$. By definition of $g_k$, the derivative $\nabla g_k(x_k)$ is surjective for all $x_k \in U$ on some neighbourhood $U$ of $\bar x$ if and only if $\nabla g(\bar x)$ is surjective; hence, $\bar x$ satisfies the necessary condition
\begin{align}
    \label{eq:stationarity}
    \nabla \varphi(\bar x) + \nabla g(\bar x)^* \bar y = 0
\end{align}
for some $\bar y$. By surjectivity, the choice of $y_k$ or $\bar y$ satisfying \eqref{eq:stationarity-quadratic} and \eqref{eq:stationarity}, respectively, is unique (for the optimal solution $\bar x_k$ or $\bar x$). Moreover, there exists a constant $\beta > 0$ such that
\begin{align}
    \label{eq:bounded-duals}
    \| y_k \| \leq \beta
\end{align}
for all $(\bar x_k, y_k) \in U \times \mathbb R^m$ satisfying \eqref{eq:stationarity-quadratic}.

Strong convexity of $\varphi_k$ relative to $\Omega_k$ is equivalent to strong monotonicity of the subgradient mapping
\begin{align*}
    \partial (\varphi_k + \imath_{\Omega_k}): x \mapsto \{ \nabla \varphi_k(x) + \nabla g_k(x)^* y \mid y \in \mathbb R^m \}
\end{align*}
where $\imath_{\Omega_k}: \mathbb R^n \to \mathbb R \cup \{\infty\}$ denotes the indicator function of the affine subspace $\Omega_k$. The necessary conditions \eqref{eq:stationarity-quadratic} can thus be written as the variational inequality $\partial (\varphi_k + \imath_{\Omega_k})(x) \ni 0$ and strong convexity of $\varphi_k$ (with constant $\mu$) implies strong regularity of the subgradient mapping (with constant $\mu^{-1}$). 

By construction of $\varphi_k$, the necessary conditions \eqref{eq:stationarity-quadratic} are a first-order approximation of \eqref{eq:stationarity}, where the gradient Lipschitz continuity condition in Assumption~\ref{ass:standing} {estimates} how well $H_k$ approximates $\nabla L(\cdot, y)$ for a given $y \in \mathbb R^m$. Monotonicity provided, convergence of sequential quadratic programming can be analyzed as generalized Newton's method for solving \eqref{eq:stationarity}. Previous works on Newton methods made assumptions of strong regularity, Lipschitz continuity of the gradient error, and a gain condition $\gamma \mu^{-1} < 1$ (equivalent to $\gamma < \mu$). 

Since the conditions of Assumption~\ref{ass:standing} are stated for the quadratic program \eqref{eq:quadratic-problem}, an algorithm may select a suitable $H_k$. Such strategies include {\em quasi}-Newton methods which approximate or regularize the Hessian of $L(\cdot, y_{k-1})$ (existence provided) around $x_k$. If $\varphi$ and $g$ are twice continuously differentiable and the optimal solution $\bar x$ satisfies a second-order sufficiency condition, Assumption~\ref{ass:standing} is satisfied by the Hessian of $L(\cdot, \bar y)$.

\begin{remark}
    The conditions in Assumption~\ref{ass:standing} could be localized to neighbourhoods of the stationary solution $(\bar x_k, y_k)$ to \eqref{eq:stationarity-quadratic}. Such a local notion of strong convexity would be equivalent to {\em local} strong monotonicity (of $\partial (\varphi_k + \imath_{\Omega_k})$) or {\em variational} strong convexity (of $\varphi_k + \imath_{\Omega_k}$) introduced in \cite{rockafellar_variational_2019}. For convex functions, both notions are equivalent to strong regularity. The results of this paper could be derived under those localized assumptions with moderate technical adaptations, which we avoid for the sake of presentation.
\end{remark}

\begin{remark}
    Surjectivity, (strong) convexity, and strong regularity of the respective mappings are also characterized by properties of the Jacobian matrix of $g$ and the standard matrix of $H_k$. We have adopted a more abstract language here to connect our work to the study of generalized Newton methods.
\end{remark}

\section{Lyapunov-like Functions}
We introduce a pair of Lyapunov-like {value functions for} optimality and constraint violation of a solution $x_k$. We then proceed to show that the solution to \eqref{eq:quadratic-problem}, which uniquely exists under Assumption~\ref{ass:standing}, jointly decreases the {value functions} in the sense of {a network dissipation inequality}. From the {dissipation inequalities} we then deduce (linear) convergence on a neighbourhood of the optimal solution $\bar x$.

{For the value of} constraint violation, we use
\begin{align*}
    v_1: x \mapsto \| g(x) \|^2
\end{align*}
and we will see that this is suitable replacement for the distance to the feasible set $\Omega$. {For the value of} optimality, we introduce
\begin{align*}
    v_2: x \mapsto \sup_{\xi \in W} \big\{ \varphi(x) - \varphi(\xi) \mid g(\xi) = g(x) \big\}
\end{align*}
for an open neighbourhood $W \subset \mathbb R^n$ of $\bar x$.
The function~$v_2$ takes the difference in the cost function relative to a {\em perturbed} optimal value. Both {value} functions are clearly nonnegative on any neighbourhood $U \subset W$ of $\bar x$. For our convergence analysis, we need to establish quadratic lower and upper bounds around $\bar x$.

\subsection{Value function for constraint violation}
Since the {value} of constraint violation is zero for any feasible solution, we can only get quadratic growth of $v_1$ relative to the feasible set $\Omega$. In the following we will repeatedly use that surjectivity of the operator $\nabla g_k(x_k): \mathbb R^n \to \mathbb R^m$ is equivalent to metric regularity of $g_k$ with some constant $\chi > 0$, or \cite[p.~38]{Dontchev2021}
\begin{align}
    \label{eq:metric-regularity-linear}
    \distance(x, g_k^{-1}(z)) \leq \chi \| z - g_k(x) \|
\end{align}
for all $x \in \mathbb R^n$ and $z \in \mathbb R^m$. By the Lyusternik--Graves theorem \cite[Theorem~5.1]{Dontchev2021} a similar, local property holds for $g$ at $\bar x$: There exists a constant $\bar\chi > 0$ along with neighbourhoods $U$ of $\bar x$ and $V$ of $g(\bar x)$ such that
\begin{align}
    \label{eq:metric-regularity}
    \distance(x, g^{-1}(z)) \leq \bar\chi \| z - g(x) \|
\end{align}
for all $x \in U$ and $z \in V$.
The quadratic bounds follow immediately from \eqref{eq:metric-regularity} and differentiability (a fortiori calmness) of $g$.

\begin{proposition}[Constraint violation]
    \label{prop:violation}
    There exist a neighbourhood $U \subset \mathbb R^n$ of $\bar x$ and constants $\zeta_1, \zeta_2 > 0$ such that
    \begin{align*}
        \zeta_1 \distance(x, \Omega)^2 \leq v_1(x) \leq \zeta_2 \| x - \bar x \|^2
    \end{align*}
    for all $x \in U$.
\end{proposition}
\begin{proof}
    By differentiability of $g$ and surjectivity of $\nabla g(\bar x)$, there exists a neighbourhood $U$ of $\bar x$ and constants $\zeta_1, \zeta_2 > 0$ such that 
    \begin{align*}
        \sqrt{\zeta_1} \, \distance(x, \Omega) \leq \| g(x) \| \leq \sqrt{\zeta_2} \, \| x - \bar x \|
    \end{align*}
    for all $x \in U$; where the first inequality follows from \eqref{eq:metric-regularity} with $z = 0$ and using that $g^{-1}(0) = \Omega$. Taking the square completes the proof.
\end{proof}

\subsection{Value function for optimality}
For the {value} of optimality, we establish a quadratic growth that is {\em complementary} to that of $v_1$ in the sense that $v_2$ grows quadratically on the feasible set $\Omega$ and behaves `sufficiently well' around $\bar x$. We first prove that the conditions of Assumption~\ref{ass:standing} imply a local notion of strong convexity also for the cost function $\varphi$ in \eqref{eq:problem} relative to its -- generally nonconvex -- feasible set. \added[id=TC]{Hence that $v_2(\bar x) = 0$, a necessary condition for optimality, is {\em equivalent} to the stationarity condition \eqref{eq:stationarity}.} 

To that extent, note that since $H_k$ is a symmetric linear operator and $\Omega_k$ is an affine subspace of $\mathbb R^n$, strong convexity of $\varphi_k$ relative to $\Omega_k$ implies also strong convexity relative to all elements of the quotient space $\mathbb R^n/\kernel \nabla g(x_k)$.

\begin{theorem}[Local strong convexity]
    \label{thm:local-convexity}
    There exist a neighbourhood $U \subset \mathbb R^n$ of $\bar x$ and a constant $\nu > 0$ such that
    \begin{align*}
        \varphi(x_1) \geq \varphi(x_2) + \langle \nabla \varphi(x_2) + \nabla g(x_2)^* \added[id=TC]{y_k}, x_1 - x_2 \rangle + \frac{\nu}{2} \| x_1 - x_2 \|^2
    \end{align*}
    for all $x_1, x_2 \in g^{-1}(b) \cap U$, \added[id=TC]{$y_k \in \Psi_k(U)$,} and $b \in \mathbb R^m$.
\end{theorem}
\begin{proof}
    Let $U \subset \mathbb R^n$ be a neighbourhood of $\bar x$ such that $\nabla [L - L_k](\cdot, \added[id=TC]{y_k})$ is Lipschitz continuous on $U$ with constant $\gamma > 0$ for all $k \in \mathbb N$ with $x_k \in U$ and \added[id=TC]{$y_k \in \Psi_k(x_k)$.} Take any $b \in \mathbb R^m$; then
    \begin{align*}
        L(x_1, \added[id=TC]{y_k}) - L(x_2, \added[id=TC]{y_k})
        \geq L_k(x_1, \added[id=TC]{y_k}) - L_k(x_2, \added[id=TC]{y_k}) + \langle \nabla [L - L_k](x_2, \added[id=TC]{y_k}), x_1 - x_2 \rangle - \frac{\gamma}{2} \| x_1 - x_2 \|^2
    \end{align*}
    for all $x_1, x_2 \in g^{-1}(b) \cap U$. Moreover, by surjectivity of $\nabla g(\bar x)$ and definition of the strict Fréchet derivative, there exist constants $\alpha, \chi > 0$ such that
    \begin{align*}
        \distance(x_1, g_k^{-1}(g_k(x_2)) \leq \chi \| g_k(x_2) - g_k(x_1) \| 
        \leq \alpha \chi \| x_1 - x_2 \|
    \end{align*}
    for all $x_1, x_2 \in U$, where the first inequality follows from \eqref{eq:metric-regularity-linear} with $z = g_k(x_2)$ and the constant $\alpha$ can be made arbitrarily small by reducing $U$. In addition, strict differentiability of \added[id=TC]{$L_k$ together with \eqref{eq:bounded-duals}} implies Lipschitz continuity of $\nabla L_k(\cdot, \added[id=TC]{y_k})$ with constant $\eta > 0$ \added[id=TC]{(independent of $y_k$)} on $U$, that is,
    \begin{align*}
        L_k(x, \added[id=TC]{y_k}) - L_k(x', \added[id=TC]{y_k}) 
        \geq \langle \nabla L_k(x', \added[id=TC]{y_k}), x - x' \rangle - \frac{\eta}{2} \| x - x' \|^2
    \end{align*}
    and
    \begin{align*}
        |\langle \nabla L_k(\hat x, \added[id=TC]{y_k}) - \nabla L_k(x', \added[id=TC]{y_k}), x - x' \rangle | \leq \eta \| x' - \hat x \| \cdot \| x - x' \|
    \end{align*}
    for all $x, x', \hat x \in U$ \added[id=TC]{and $y_k \in \Psi_k(U)$}.

    Combining above inequalities with strong convexity of $\varphi_k$ relative to $g_k^{-1}(g_k(x_2))$, for any $x_1, x_2 \in g^{-1}(b) \cap U$, there exists a solution $x_1' \in g_k^{-1}(g_k(x_2))$ satisfying
    \begin{align*}
        \| x_1' - x_1 \| \leq \alpha \chi \| x_1 - x_2 \|
    \end{align*}
    and hence, assuming without loss of generality that $x_1' \in U$ and using the triangle inequalities
    \begin{align*}
        \| x_1 - x_2 \| \leq \| x_1' - x_2 \| + \| x_1' - x_1 \| \\
        \| x_1' - x_2 \| \leq \| x_1  -x_2 \| + \| x_1' - x_1 \|
    \end{align*}
    we have that
    \begin{multline*}
        \varphi(x_1) - \varphi(x_2) - \langle \nabla L(x_2, \added[id=TC]{y_k}), x_1 - x_2 \rangle \\
    \begin{aligned}
        &\geq L_k(x_1', \added[id=TC]{y_k}) - L_k(x_2, \added[id=TC]{y_k}) - \langle \nabla L_k(x_2, \added[id=TC]{y_k}), x_1 - x_2 \rangle
        + L_k(x_1, \added[id=TC]{y_k}) - L_k(x_1', \added[id=TC]{y_k}) - \frac{\gamma}{2} \| x_1 - x_2 \|^2 \\
        &\begin{multlined}[c][.8\linewidth]
            = L_k(x_1', \added[id=TC]{y_k}) - L_k(x_2, \added[id=TC]{y_k}) - \langle \nabla L_k(x_2, \added[id=TC]{y_k}), x_1' - x_2 \rangle \\
            + L_k(x_1, \added[id=TC]{y_k}) - L_k(x_1', \added[id=TC]{y_k}) + \langle \nabla L_k(x_2, \added[id=TC]{y_k}), x_1' - x_1 \rangle 
            - \frac{\gamma}{2} \| x_1 - x_2 \|^2 
        \end{multlined} \\
        &\geq \frac{\mu}{2} \| x_1' - x_2 \|^2 - \frac{\eta}{2} \| x_1 - x_1' \|^2 - \eta \| x_1' - x_2 \| \cdot \| x_1 - x_1' \| - \frac{\gamma}{2} \| x_1 - x_2 \|^2 \\
        &\begin{multlined}[c][.8\linewidth]
            \geq \frac{\mu}{2} \big( \| x_1 - x_2 \| - \alpha \chi \| x_1 - x_2 \| \big)^2 - \frac{\eta}{2} (\alpha \chi)^2 \| x_1 - x_2 \|^2 \\
            - \eta \big( \| x_1 - x_2 \| + \alpha \chi \| x_1 - x_2 \| \big) \alpha \chi \| x_1 - x_2 \| - \frac{\gamma}{2} \| x_1 - x_2 \|^2 
        \end{multlined}
    \end{aligned} \\
        \geq \frac{\nu}{2} \| x_1 - x_2 \|^2
    \end{multline*}
    where $\alpha$ can be reduced by shrinking $U$ such that
    \begin{align*}
        \frac{\nu}{2} := \frac{\mu}{2} (1 - \alpha \chi)^2 - \frac{\gamma}{2} - \eta(\alpha \chi)(1 + \alpha \chi) - \frac{\eta}{2}(\alpha \chi)^2 > 0
    \end{align*}
    since $\mu - \gamma > 0$. As the constants $\alpha, \eta, \mu, \chi$, and hence $\nu$, are independent of the choice of $b$, the proof is complete.
\end{proof}

Note that, by surjectivity of $\nabla g(\bar x)$, for any neighbourhood $U \subset \mathbb R^n$ of $\bar x$ there exists a neighbourhood $V \subset \mathbb R^m$ of $0$ such that $g^{-1}(\cdot) \cap U$ is nonempty on $V$ \cite[Proposition~4.1]{Dontchev2021}.
Since $\Omega = g^{-1}(0)$ \added[id=TC]{and $\bar y \in \Psi_k(\bar x)$}, Theorem~\ref{thm:local-convexity} implies, amongst others, that $\varphi$ grows quadratically on $\Omega$ relative to and around $\bar x$. 
The quadratic bounds on $v_2$ are further consequences of the local strong convexity notion.

\begin{proposition}[Optimality]
    \label{prop:optimality}
    There exist neighbourhoods $U, W \subset \mathbb R^n$ of $\bar x$ and constants $\zeta_3, \zeta_4, \varrho > 0$ such that
    \begin{align*}
        \zeta_3 \| x - \bar x \|^2 - \varrho \distance(x, \Omega)^2
        \leq v_2(x) \leq \zeta_4 \| x - \bar x \|^2
    \end{align*}
    and $v_2(x) \geq 0$ for all $x \in U$.
\end{proposition}
\begin{proof}
    By virtue of Theorem~\ref{thm:local-convexity}, there exists a neighbourhood $W \subset \mathbb R^n$ such that $\varphi$ satisfies a strong convexity condition on $g^{-1}(b) \cap W$ with some constant $\nu > 0$ for all $b \in \mathbb R^m$.
    Moreover, strict differentiability of $L(\cdot, \bar y)$ implies Lipschitz continuity of $\nabla L(\cdot, \bar y)$ with constant $\bar \eta > 0$ on $W$ and $g$ is, a fortiori, Lipschitz continuous with constant $\bar \lambda > 0$ on $W$.
    Choose a neighbourhood $U \subset W$ of $\bar x$.

    For the lower bound, take $x \in U$; by surjectivity of $\nabla g(\bar x)$, there exists a solution $\xi \in g^{-1}(g(x))$ satisfying
    \begin{align*}
        \| \xi - \bar x \| \leq \bar\chi \| g(x) \|
    \end{align*}
    which is \eqref{eq:metric-regularity} with $z = g(x)$, $g(\bar x) = 0$, and $\xi \in W$ without loss of generality. Hence,
    \begin{align*}
        v_2(x) \geq \varphi(x) - \varphi(\xi) \geq \langle \nabla L(\xi, \bar y), x - \xi \rangle + \frac{\nu}{2} \| x - \xi \|^2 
        \geq - \bar \eta \| \xi - \bar x \| \cdot \| x - \xi \| + \frac{\nu}{2} \| x - \xi \|^2
    \end{align*}
    where the first inequality holds since $\xi$ is a feasible solution. Using Lipschitz continuity of $g$ and the triangle inequalities
    \begin{align*}
        \| x - \xi \| \leq \| \xi - \bar x \| + \| x - \bar x \| \\
        \| x - \bar x \| \leq \| \xi - \bar x \| + \| x - \xi \|
    \end{align*}
    {there exist constants $\zeta_3, \hat \varrho > 0$ such that}
    the lower bound can be rewritten to\footnotemark
    \begin{multline*}
        v_2(x) 
        \geq \frac{\nu}{2} \| x - \bar x \|^2 - (\bar \eta + \nu) \| \xi - \bar x \| \cdot \| x - \bar x \| - (\bar \eta - \frac{\nu}{2}) \| \xi - \bar x \|^2 \\
        \geq \zeta_3 \| x - \bar x \|^2 - \hat \varrho \| \xi - \bar x \|^2
        \geq \zeta_3 \| x - \bar x \|^2 - \varrho \distance(x, \Omega)^2
    \end{multline*}
    with $\varrho := \added[id=IK]{(\bar\chi \bar \lambda)^2 \hat \varrho}$. 
    \footnotetext{Specifically, $2\zeta_3 = \nu + \varepsilon (\bar\eta + \nu)$ and $2\hat\varrho = 2\bar\eta + \nu + \varepsilon^{-1} (\bar\eta + \nu)$ with $\varepsilon > 0$ by Young's inequality.}

    For the upper bound, 
    take $x \in U$ and $\xi \in g^{-1}(g(x)) \cap W$; by surjectivity of $\nabla g(\bar x)$, there exists a solution $\xi' \in g^{-1}(0)$ satisfying
    \begin{align*}
        \| \xi' - \xi \| \leq \bar\chi \| g(\xi) \| = \bar\chi \| g(x) \|
    \end{align*}
    which is \eqref{eq:metric-regularity} with $z = 0$, and $\xi' \in W$ without loss of generality. Using Lipschitz continuity of $g$ and the triangle inequalities
    \begin{align*}
        \| \xi - \bar x \| \leq \| \xi' - \bar x \| + \| \xi' - \xi \| \\
        \| \xi' - \bar x \| \leq \| \xi - \bar x \| + \| \xi' - \xi \|
    \end{align*}
    {there exists a constant $\omega > 0$ such that}
    we arrive at\footnotemark
    \begin{multline*}
        L(\xi, \bar y) \geq \varphi(\bar x) - \langle \nabla L(\xi, \bar y), \xi' - \xi \rangle
        + \frac{\nu}{2} \| \xi' - \bar x \|^2 - \frac{\eta}{2} \| \xi' - \xi \|^2 \\
        \geq \varphi(\bar x) + \frac{\nu}{2} \big( \| \xi - \bar x \| - \bar\chi \| g(x) \| \big)^2  
        - \bar \eta \| \xi - \bar x \| \cdot \bar\chi \| g(x) \| - \frac{\bar \eta}{2} \bar\chi^2 \| g(x) \|^2 \\
        \geq \varphi(\bar x) - \frac{\omega}{2} \| g(x) \|^2
        \geq \varphi(\bar x) - \frac{\bar \omega}{2} \| x - \bar x \|^2
    \end{multline*}
    with $\bar \omega := \bar \lambda^2 \omega$. 
    \footnotetext{Specifically, $\omega = (\bar\eta - \nu) \bar\chi^2 + \varepsilon (\nu + \bar\eta) \bar\chi$ with $\varepsilon > 0$ by Young's inequality, assuming w.l.o.g. that $\bar\eta > \nu$.}
    On the other hand, Lipschitz continuity of $\nabla L(\cdot, \bar y)$ gives
    \begin{align*}
        L(x, \bar y) \leq L(\bar x, \bar y) + \langle \nabla L(\bar x, \bar y), x - \bar x \rangle + \frac{\bar \eta}{2} \| x - \bar x \|^2
    \end{align*}
    for all $x \in U$; and note that $\nabla L(\bar x, \bar y) = 0$. Combining these inequalities, we have that
    \begin{align*}
        \varphi(x) - \varphi(\xi) = L(x, \bar y) - L(\xi, \bar y) 
        \leq \big(\frac{\bar \eta}{2} + \frac{\bar \omega}{2}\big) \| x - \bar x \|^2
    \end{align*}
    for all $x \in U$ and $\xi \in g^{-1}(g(x)) \cap W$. Taking the supremum over $\xi \in W$ provides the desired upper bound with $\zeta_4 = \frac{\bar \eta + \bar \omega}{2}$.

    Since $U \subset W$, there exist a feasible solution $\xi \in W$ subject to $g(\xi) = g(x)$ for any $x \in U$; thus $v_2(x) \geq 0$ and the proof is complete.
\end{proof}

\section{Dissipation Inequalities}
We now show that the update $x_{k+1} \in \Phi_k(x_k)$ satisfies {ISS}-type dissipation inequalities for the {value functions} for constraint violation and optimality. In particular, we show that passing from $x_k$ to $x_{k+1}$ decreases $v_1$ \added[id=TC]{and $v_2$ if $x_k$ is already close to optimality or the feasible set, respectively; otherwise, their values might increase, but always remain bounded}. This result, similar to dissipation arguments for networks of subsystems \cite{Arcak2016}, is the cornerstone of deriving a strict Lyapunov function that establishes asymptotic stability and linear convergence in the next section.

We start by noting a Polyak--Łojasiewicz (PL)-type inequality and bound on the distance between solutions to the perturbed optimization problem as immediate consequences of Theorem~\ref{thm:local-convexity}.

\begin{corollary}[PL-type inequality]
    \label{cor:pl-inequality}
    Take any $b \in \mathbb R^m$, \added[id=TC]{$y_k \in \Psi_k(U)$,} and $x_1, x_2 \in g^{-1}(b) \cap U$; then
    \begin{align*}
        \varphi(x_2) - \varphi(x_1) \leq \frac{1}{2\nu} \| \nabla \varphi(x_2) + \nabla g(x_2)^* \added[id=TC]{y_k} \|^2
    \end{align*}
    where $U$ and $\nu$ are defined in Theorem~\ref{thm:local-convexity}.
\end{corollary}
\begin{proof}
    We have that
    \begin{align*}
        \varphi(x_1) \geq \varphi(x_2) - \| \nabla \varphi(x_2) + \nabla g(x_2)^* \added[id=TC]{y_k} \| \cdot \| x_1 - x_2 \| + \frac{\nu}{2} \| x_1 - x_2 \|^2
    \end{align*}
    by virtue of Theorem~\ref{thm:local-convexity}. Minimizing the right-hand side with respect to $\| x_1 - x_2 \|$ yields
    \begin{align*}
        \varphi(x_1) - \varphi(x_2) \geq -\frac{1}{2\nu} \| \nabla \varphi(x_2) + \nabla g(x_2)^* \added[id=TC]{y_k} \|^2
    \end{align*}
    and negation gives the desired result.
\end{proof}

\begin{corollary}[Distance bound]
    \label{cor:distance-bound}
    Take any $b \in \mathbb R^m$, \added[id=TC]{$y_k \in \Psi_k(U)$,} and $x_1, x_2 \in g^{-1}(b) \cap U$ \added[id=TC]{with} $\varphi(x_1) \leq \varphi(x_2)$; then
    \begin{align*}
        \| x_1 - x_2 \| \leq \frac{2}{\nu} \| \nabla \varphi(x_2) + \nabla g(x_2)^* \added[id=TC]{y_k} \|
    \end{align*}
    where $U$ and $\nu$ are defined in Theorem~\ref{thm:local-convexity}.
\end{corollary}
\begin{proof}
    Suppose that, \added[id=TC]{on the contrary,} 
    \begin{align*}
        \| x_1 - x_2 \| > \frac2\nu \| \nabla L(x_2, y_k) \|;
    \end{align*}
    then
    \begin{align*}
        \varphi(x_1) \geq \varphi(x_2) - \| \nabla L(x_2, \added[id=TC]{y_k}) \| \cdot \| x_1 - x_2 \| + \frac{\nu}{2} \| x_1 - x_2 \|^2 
        > \varphi(x_2)
    \end{align*}
    by virtue of Theorem~\ref{thm:local-convexity} \added[id=TC]{and the proof is complete.}
\end{proof}

We are ready to establish our main result, the dissipation inequalities for $v_1$ and $v_2$.

\begin{proposition}[Dissipation inequalities]
    \label{prop:matrosov}
    {There exist a neighbourhood $W \subset \mathbb R^n$ and constants $\vartheta_1, \vartheta_2, \psi_2 > 0$ such that for all $\varkappa > 0$, there exists a neighbourhood $U \subset W$ and a constant $\psi_1 > 0$ satisfying}
    \begin{align*}
        v_1(x_{k+1}) - v_1(x_k) &\leq -\vartheta_1 v_1(x_k) + \psi_1 v_2(x_k) \\
        v_2(x_{k+1}) - v_2(x_k) &\leq -\vartheta_2 v_2(x_k) + \psi_2 v_1(x_k)
    \end{align*}
    and
    \begin{align*}
        \vartheta_2 - \varkappa \psi_1 > 0
    \end{align*}
    for all $x_k \in U$ and $x_{k+1} \in \Phi(x_k) \cap U$.
\end{proposition}
\begin{proof}
    Pick a neighbourhood $W \subset \mathbb R^n$ of $\bar x$ such that the gradient Lipschitz continuity condition of Assumption~\ref{ass:standing} and the assertion of Theorem~\ref{thm:local-convexity} are satisfied over $W$.
    Choose a neighbourhood $U \subset W$ of $\bar x$.
    
    Take $x_k \in U$ and $x_{k+1} = \bar x_k$, where $(\bar x_k, \bar y_k) \in (\Phi_k(x_k) \cap U) \times \added[id=TC]{\Psi_k(x_k)}$ is a primal-dual solution of \eqref{eq:quadratic-problem}. Observe that
    \begin{align*}
        g_k(x_k) = g(x_k), \quad\quad g_k(\bar x_k) = 0
    \end{align*}
    by definition of \eqref{eq:quadratic-problem} and hence, there exists a constant $\alpha > 0$ such that
    \begin{align*}
        \| g(\bar x_k) \| \leq \alpha \| x_k - \bar x_k \|
    \end{align*}
    and hence,
    \begin{align}
        \label{eq:matrosov-violation}
        \| g(\bar x_k) \|^2 - \| g(x_k) \|^2 \leq \alpha^2 \| x_k - \bar x_k \|^2 - \| g(x_k) \|^2 
    \end{align}
    where the constant $\alpha$ can be made arbitrarily small by reducing $U$. 
    %
    %
    Moreover, following the arguments in the proof of Proposition~\ref{prop:violation}, there exist constants $\hat\zeta, \hat\varrho, \bar\chi > 0$ such that
    \begin{align*}
        v_2(x) \geq \hat\zeta \| x - \bar x \|^2 - \hat\varrho \bar\chi^2 \| g(x) \|^2
    \end{align*}
    for all $x \in U$. 
    \added[id=TC]{Strong convexity of $\varphi_k$ relative to $\Omega_k$ gives (see Appendix)
    \begin{align*}
        \| x_k - \bar x_k \| \leq \frac{1}{\mu} \| \nabla L_k(x_k, \bar y) \| + \big( \frac{\eta}{\mu} + 1 \big) \chi \| g_k(x_k) \| 
        \leq \frac{\bar\eta}{\mu} \| x_k - \bar x \| + \big( \frac{\eta}{\mu} + 1 \big) \chi \| g(x_k) \|
    \end{align*}
    where the second inequality follows from $\nabla L_k(x_k, \bar y) = \nabla L(x_k, \bar y)$ and $g_k(x_k) = g(x_k)$, Lipschitz continuity of $\nabla L(\cdot, \bar y)$ with constant $\bar\eta$, and $\nabla L(\bar x, \bar y) = 0$.}
    Using \eqref{eq:matrosov-violation} \added[id=TC]{and} the triangle inequality
    \begin{align*}
        \| x_k - \bar x_k \| \leq \| x_k' - \bar x_k \| + \| x_k' - x_k \|,
    \end{align*}
    there exist constants $\hat\psi, \omega > 0$ such that\footnotemark
    \begin{multline*}
        v_1(x_{k+1}) - v_1(x_k) \\
    \begin{aligned}
        &\leq \alpha^2 \big(\frac{\bar\eta}{\mu}\big)^2 \| x_k - \bar x \|^2
        + \added[id=TC]{2 \alpha^2 \frac{\bar\eta}{\mu} \big( \frac{\eta}{\mu} + 1 \big) \chi \| x_k - \bar x \| \cdot \| g(x_k) \|}
        + \added[id=TC]{\alpha^2 \big( \frac{\eta}{\mu} + 1 \big)^2 \chi^2 \| g(x_k) \|^2 - \| g(x_k) \|^2} \\
        &\leq \added[id=TC]{\alpha^2 \hat\zeta \hat\psi \| x_k - \bar x \|^2 + (\alpha^2 \omega \chi^2 - 1) \| g(x_k) \|^2}
    \end{aligned} \\
        \leq \alpha^2 \hat\psi v_2(x_k) + (\alpha^2 \hat\varrho \bar\chi^2 \hat\psi + \alpha^2 \chi^2 \omega - 1) \| g(x_k) \|^2
        = -\vartheta_1 v_1(x_k) + \psi_1 v_2(x_k)
    \end{multline*}
    with $\psi_1 := \alpha^2 \hat\psi$, where
    \footnotetext{Specifically, 
    \[
        \hat\zeta \hat\psi = \big( \frac{\bar\eta}{\mu} \big)^2 + \varepsilon \frac{\bar\eta}{\mu} \big( 1 + \frac{\eta}{\mu} \big), \qquad \omega = \big( 1 + \frac{\eta}{\mu} \big)^2 + \varepsilon^{-1} \frac{\bar\eta}{\mu} \big( 1 + \frac{\eta}{\mu} \big)
    \]
    with $\varepsilon = \chi \varepsilon' > 0$ by Young's inequality.}
    $\alpha$ can be reduced by shrinking $U$ such that $\vartheta_1 := 1 - \alpha^2 (\hat\varrho \bar\chi^2 \hat\psi + \chi^2 \omega) > \added[id=IK]{0}$.


    By virtue of Corollary~\ref{cor:pl-inequality}, we have that
    \begin{align*}
        \varphi(x_k) - \varphi(\xi_k) \leq \frac{1}{2\nu} \| \nabla L_k(x_k, \bar y_k) \|^2
    \end{align*}
    for all $\xi_k \in g^{-1}(g(x_k)) \cap W$, using that $\nabla L(x_k, \bar y_k) = \nabla L_k(x_k, \bar y_k)$.
    Taking the supremum 
    over $\xi_k$, we obtain
    \begin{align*}
        v_2(x_k) \leq \frac{\eta^2}{2 \nu} \| x_k - \bar x_k \|^2 
    \end{align*}
    where we use Lipschitz continuity of $\nabla L_k(\cdot, \bar y_k)$ with constant $\eta$, and $\nabla L_k(\bar x_k, \bar y_k) = 0$.
    
    Take any $\xi \in g^{-1}(g(\bar x_k)) \cap W$ \added[id=TC]{with $\varphi(\xi) \leq \varphi(\bar x_k)$}.
    \added[id=TC]{By virtue of Corollary~\ref{cor:distance-bound}, we have that
    \begin{align*}
        \| \xi - \bar x_k \| < \frac{2}{\nu} \| \nabla L(\bar x_k, \bar y_k) \| \leq \frac{2\gamma}{\nu} \| x_k - \bar x_k \|
    \end{align*}
    where the second inequality follows from Lipschitz continuity of $\nabla [L - L_k]$ with constant $\gamma$.
    }
    By surjectivity of $\nabla g(\bar x)$, there exists a solution $\xi' \in g^{-1}(g(x_k))$ satisfying
    \begin{align*}
        \| \xi' - \xi \| \leq \bar\chi \| g(\bar x_k) - g(x_k) \|
    \end{align*}
    which is \eqref{eq:metric-regularity} with $z = g(x_k)$; and \added[id=TC]{by shrinking $U$ we can ensure that $\xi' \in W$}.
    Moreover, there exists a solution $x_k' \in g_k^{-1}(0)$ satisfying
    \begin{align*}
        \| x_k' - x_k \| \leq \chi \| g_k(x_k) \| = \chi \| g(x_k) \|
    \end{align*}
    which is \eqref{eq:metric-regularity-linear} with $z = 0$.
    Since $\varphi_k$ is strongly convex on $\Omega = g_k^{-1}(0)$ with constant $\mu > 0$ and $\nabla L_k(\cdot, \bar y_k)$ is Lipschitz continuous with constant $\eta > 0$, as well as using the triangle inequalities
    \begin{align*}
        \| x_k - \bar x_k \| &\leq \| x_k' - \bar x_k \| + \| x_k' - x_k \| \\
        \| \xi - x_k \| &\leq \| x_k - \bar x_k \| + \| \xi - \bar x_k \| \\
        \| g(\bar x_k) - g(x_k) \| &\leq \| g(\bar x_k) \| + \| g(x_k) \|
    \end{align*}
    and
    \begin{align*}
        \| \xi' - \xi \| &\leq \alpha \bar\chi \| x_k - \bar x_k \| + \bar\chi \| g(x_k) \| \\
        \| \xi - \bar x_k \| 
            &\leq \frac{2\gamma}{\nu} \| x_k' - \bar x_k \| + \frac{2\gamma}{\nu} \chi \| g(x_k) \|
    \end{align*}
    we have that 
    \begin{multline}
        \label{eq:matrosov-optimality}
        \varphi(\bar x_k) - \varphi(\xi) - [\varphi(x_k) - \varphi(\xi')]
        = L(\bar x_k, \bar y_k) - L(x_k, \bar y_k) + L(\xi', \bar y_k) - L(\xi, \bar y_k) \\
        \begin{multlined}[c][.85\linewidth]
            \leq L_k(\bar x_k, \bar y_k) - L_k(x_k, \bar y_k) \added[id=IK]{+} \langle \nabla [L - L_k](x_k, \bar y_k), \bar x_k - x_k \rangle + \frac{\gamma}{2} \| x_k - \bar x_k \|^2 \\
            + L_k(\xi', \bar y_k) - L_k(\xi, \bar y_k) \added[id=IK]{+} \langle \nabla [L - L_k](\xi, \bar y_k), \xi' - \xi \rangle + \frac{\gamma}{2} \| \xi' - \xi \|^2 
        \end{multlined} \\
        \begin{multlined}[c][.85\linewidth]
            \leq \frac{\gamma}{2} \| x_k - \bar x_k \|^2 - \frac{\mu}{2} \| x_k' - \bar x_k \|^2 + \frac{\eta}{2} \| x_k' - x_k \|^2 + \eta \| x_k' - x_k \| \cdot \| x_k' - \bar x_k \| \\
            + \frac{\gamma + \eta}{2} \| \xi' - \xi \|^2 + (\gamma + \eta) \| \xi' - \xi \|    \cdot \| \xi - x_k \| 
        \end{multlined} \\
        \begin{multlined}[c][.85\linewidth]
            \leq \frac{\gamma - \mu}{2} \| x_k' - \bar x_k \|^2 + (\gamma + \eta) \| x_k' - \bar x_k \| \cdot \| x_k' - x_k \| + \frac{\gamma + \eta}{2} \| x_k' - x_k \|^2 \\
            + \frac{\gamma + \eta}{2} \| \xi' - \xi \|^2 + (\gamma + \eta) \| \xi' - \xi \| \cdot \| x_k - \bar x_k \| + (\gamma + \eta) \| \xi' - \xi \| \cdot \| \xi - \bar x_k \|
        \end{multlined} \\
        \begin{multlined}[c][.85\linewidth]
            \leq \frac{\gamma - \mu}{2} \| x_k' - \bar x_k \|^2 + (\gamma + \eta) \chi \| x_k' - \bar x_k \| \cdot \| g(x_k) \| + \frac{\gamma + \eta}{2} \chi^2 \| g(x_k) \|^2 \\
        \begin{aligned}
            &+ \frac{\gamma + \eta}{2} \Big[ (\alpha \bar\chi)^2 \| x_k' - \bar x_k \|^2 + 2 \alpha \bar\chi (\bar\chi + \chi) \| x_k' - \bar x_k \| \cdot \| g(x_k) \| + (\bar\chi + \chi)^2 \| g(x_k) \|^2 \Big] \\
            &+ (\gamma + \eta) \Big[ \alpha \bar\chi \| x_k' - \bar x_k \|^2 + (2 \alpha \chi + 1) \bar\chi \| x_k' - \bar x_k \| \cdot \| g(x_k) \| + \chi \bar\chi \| g(x_k) \|^2 \Big] \\
            &+ (\gamma + \eta) \Big[ \alpha \bar\chi \frac{2\gamma}{\nu} \| x_k' - \bar x_k \|^2 + (\alpha \chi \bar\chi + \bar\chi + \chi) \frac{2\gamma}{\nu} \| x_k' - \bar x_k \| \cdot \| g(x_k) \| + \frac{2\gamma}{\nu} \chi (\bar\chi + \chi) \| g(x_k) \|^2 \Big]
        \end{aligned}
        \end{multlined} \\
        %
        %
        = -\pi \| x_k' - \bar x_k \|^2 + 2 \rho \| x_k' - \bar x_k \| \cdot \| g(x_k) \| + \sigma \| g(x_k) \|^2
    \end{multline}
    with constants 
    \begin{align*}
        \pi &
        = \frac{\mu - \gamma}{2} - (\gamma + \eta) \bar\chi \alpha \Big[\frac{1}{2} \bar\chi \alpha + \big(1 + \frac{2\gamma}{\nu}\big) \Big] \\
        \rho &= \frac{1}{2} (\gamma + \eta) \Big[ \chi + \alpha \bar\chi (\bar\chi + \chi) + (2 \alpha \chi + 1) \bar\chi + (\alpha \chi \bar\chi + \bar\chi + \chi) \frac{2 \gamma}{\nu} \Big] \\
        \sigma &= (\gamma + \eta) \Big[ \frac{1}{2} \chi^2 + \frac{1}{2} (\bar\chi + \chi)^2 + \chi \bar\chi + \frac{2 \gamma}{\nu} \chi (\bar\chi + \chi) \Big]
    \end{align*}
    where $\rho, \sigma > 0$ and, since $\mu > \gamma$,
    we can reduce $\alpha$ by shrinking $U$ to ensure that $\pi > 0$.

    Taking the supremum of \eqref{eq:matrosov-optimality} over $\xi$, noting that $v_2(x_k) \geq \varphi(x_k) - \varphi(\xi')$, there hence exist constants $\hat\vartheta, \psi_2 > 0$ such that\footnotemark
    \begin{multline*}
        v_2(x_{k+1}) - v_2(x_k) \\ 
        \leq -\pi \big( \| x_k' - \bar x_k \| + \chi \| g(x_k) \| \big)^2
        + 2 (\rho + \pi \chi) \| x_k' - \bar x_k \| \cdot \| g(x_k) \| + (\sigma + \pi \chi^2) \| g(x_k) \|^2 \\
        \leq -\hat\vartheta \big( \| x_k' - \bar x_k \| + \chi \| g(x_k) \| \big)^2 + \psi_2 \| g(x_k) \|^2 
        \leq -\vartheta_2 v_2(x_k) + \psi_2 v_1(x_k)
    \end{multline*}
    with $\vartheta_2 := \hat\vartheta \frac{2 \nu}{\eta^2}$. 
    \footnotetext{Specifically,
    \[
        \hat\vartheta = \pi - \frac{\varepsilon}{1 + \varepsilon \chi} (\rho + \pi \chi), \quad
        \psi_2 = (\sigma + \pi \chi^2) + (\varepsilon^{-1} - \chi) (\rho + \pi \chi)
    \]
    with $\varepsilon \in (0, \chi^{-1})$ by Young's inequality.
    }

    Note that the constants $\vartheta_1, \vartheta_2, \hat\psi, \psi_2$ chosen so far satisfy their claimed properties if $\alpha$ is further reduced. Pick now any $\varkappa > 0$ and reduce $\alpha$ (and hence $\psi_1 = \alpha^2 \hat\psi$) by shrinking $U$ such that $\vartheta_2 - \varkappa \psi_1 > 0$. This completes the proof.
\end{proof}

\section{Stability \& Convergence Analysis}
We now establish stability and linear convergence of the Newton iteration
\begin{align}
    x_{k+1} \in \Phi_k(x_k), \quad k \in \mathbb N, \tag{\ref{eq:newton-iteration}}
\end{align}
in an neighbourhood \added[id=TC]{$U$} of $\bar x$.
\added[id=TC]{To that extent, we first construct} a strict local Lyapunov function with quadratic bounds \added[id=TC]{and derive linear convergence on $U$. Our final theorem then completes the analysis by certifying that the sequence generated by \eqref{eq:newton-iteration} remains in $U$ if started sufficiently close to $\bar x$.}

\begin{lemma}[Strict Lyapunov function]
    \label{lem:strict-lyapunov}
    There exist neighbourhoods $U, W \subset \mathbb R^n$ of $\bar x$ and constants $\varkappa, c_0, c_1, c_2 > 0$ such that the function
    \begin{align*}
        S = \varkappa v_1 + v_2
    \end{align*}
    satisfies
    \begin{align*}
        c_1 \| x - \bar x \|^2 \leq S(x) &\leq c_2 \| x - \bar x \|^2 \\
        S(x_{k+1}) - S(x_k) &\leq -c_0 S(x_k)
    \end{align*}
    for all $x, x_k \in U$ and $x_{k+1} \in \Phi(x_k) \cap U$.
\end{lemma}
\begin{proof}
    Pick a neighbourhood $W \subset \mathbb R^n$ of $\bar x$ such that the assertions of Propositions~\ref{prop:violation}--\ref{prop:matrosov} with constants $\zeta_1, \zeta_2, \zeta_3, \zeta_4, \varrho > 0$ and $\vartheta_1, \vartheta_2, \psi_2 > 0$ are satisfied over some (joint) neighbourhood $\bar U \subset W$.
    Take $\varkappa > 0$ and $c_1, c_2 > 0$ such that
    \begin{align*}
        \varkappa > \max\{ \tfrac{\varrho}{\zeta_1}, \tfrac{\psi_2}{\vartheta_1} \}, \quad c_1 \leq \min \{ \varkappa \tfrac{\zeta_1 \zeta_3}{\varrho}, \zeta_3 \},
    \end{align*}
    and $c_2 \geq \zeta_2 \varkappa + \zeta_4$. Then, by virtue of Propositions~\ref{prop:violation} and \ref{prop:optimality},
    \begin{align*}
        S(x) \geq (\varkappa \zeta_1 - \varrho) \distance(x, \Omega)^2 + \zeta_3 \| x - \bar x \|^2
    \end{align*}
    \added[id=IK]{since $\varkappa \geq \tfrac{\varrho}{\zeta_1}$. Hence, with $c_1 \leq \zeta_3$,}
    \begin{align*}
        c_1 \| x - \bar x \|^2 \leq S(x) \leq (\varkappa \zeta_2 + \zeta_4) \| x - \bar x \|^2 \leq c_2 \| x - \bar x \|^2
    \end{align*}
    for all $x \in \bar U$. 
    Moreover, by virtue of Proposition~\ref{prop:matrosov}, there exists a neighbourhood $U' \subset W$ and constant $\psi_1 > 0$ such that
    \begin{multline*}
        S(x_{k+1}) - S(x_k)
        \leq \added[id=IK]{-\varkappa \big[ \vartheta_1 v_1(x_k) - \psi_1 v_2(x_k) \big] - \big[ \vartheta_2 v_2(x_k) - \psi_2 v_1(x_k) \big]} \\
        = -(\vartheta_1 - \frac{\psi_2}{\varkappa}) \varkappa v_1(x_k) - (\vartheta_2 - \varkappa \psi_1) v_2(x_k)
        \leq -c_0 \big[ \varkappa v_1(x_k) + v_2(x_k) \big]
        = -c_0 S(x_k)
    \end{multline*}
    with constant $c_0 > 0$ satisfying
    \begin{align*}
        c_0 \leq \min\{ (\vartheta_1 - \tfrac{\psi_2}{\varkappa}), (\vartheta_2 - \varkappa \psi_1) \}
    \end{align*}
    for all $x_k \in U'$ and $x_{k+1} \in \Phi(x_k) \cap U'$. Choosing a neighbourhood $U \subset \bar U \cap U'$ completes the proof.
\end{proof}

\added[id=TC]{Since the Lyapunov function $S(\cdot)$ has quadratic bounds, we show linear convergence (in the root sense).}

\begin{corollary}[R-Linear convergence]
    \label{cor:linear-convergence}
    There exist a neighbourhood $U \subset \mathbb R^n$ of $\bar x$ and constants $\Gamma > 0$, $p \in (0, 1)$ 
    such that
    \begin{align*}
        \forall k \in \mathbb N, \qquad \| x_k - \bar x \| \leq \Gamma p^k \| x_0 - \bar x \|
    \end{align*}
    for any sequence $(x_k)_{k \in \mathbb N} \subset U$ generated by \eqref{eq:newton-iteration}.
\end{corollary}
\begin{proof}
    Let the assertions of Lemma~\ref{lem:strict-lyapunov} be satisfied by the function $S: \mathbb R^n \to \mathbb R$ with constants $c_0, c_1, c_2 > 0$ on an neighbourhood $U \subset \mathbb R^n$ of $\bar x$. Any sequence $(x_k)_{k \in \mathbb N} \subset U$ generated by \eqref{eq:newton-iteration} satisfies
    \begin{align*}
        c_1 \| x_k - \bar x \|^2 \leq 
        S(x_k)
        \leq (1 - c_0) S(x_{k-1}) \leq \cdots \leq (1 - c_0)^k S(x_0)
        \leq (1 - c_0) c_2 \| x_0 - \bar x \|^2
    \end{align*}
    for all $k \in \mathbb N$. Rewriting with $\Gamma = \sqrt{c_2/c_1}$ and $p = \sqrt{1 - c_0}$ completes the proof.
\end{proof}

\added[id=TC]{Finally, we establish invariance of a neighbourhood of $\bar x$ and thus prove local exponential stability.}

\begin{theorem}[Asymptotic stability]
    For any neighbourhood $U_1 \subset \mathbb R^n$ of $\bar x$, there exists a neighbourhood $U_0 \subset U_1$ of $\bar x$ such that the sequence $(x_k)_{k \in \mathbb N}$ generated by \eqref{eq:newton-iteration} for $x_0 \in U_0$ (uniquely) exists, satisfies $x_k \in U_1$ for all $k \in \mathbb N$, and converges linearly to $\bar x$.
\end{theorem}
\begin{proof}
    Let the assertions of Corollary~\ref{cor:linear-convergence} be satisfied with constants $\Gamma > 0$ and $p \in (0, 1)$ on an neighbourhood $U \subset \mathbb R^n$ of $\bar x$. By virtue of the standing assumptions, there exists a (unique) primal-dual solution $(\bar x_k, \bar y_k) \in \mathbb R^n \times \mathbb R^m$ of \eqref{eq:quadratic-problem} for any $x_k \in U$; 
    and $[g - g_k]$ is Lipschitz continuous with constant $\lambda > 0$ on $U$. Since $\varphi_k$ is strongly convex relative to $\Omega_k$, we have that \added[id=TC]{(see Appendix)}
    \begin{align*}
        \| \bar x_k - \bar x \| \leq \frac{1}{\mu} \| \nabla L_k(\bar x, \bar y) \| + \big( \frac{\eta}{\mu} + 1 \big) \chi \| g_k(\bar x) \|
        \leq \big( \frac{\gamma}{\mu} + \frac{\eta}{\mu} \lambda \chi + \lambda \chi \big) \| x_k - \bar x \| =: \hat\delta \| x_k - \bar x \|
    \end{align*}
    where we have used that $\nabla L(\bar x, \bar y) = 0$ and $g(\bar x) = 0$.
    Given a neighbourhood $U_1 \subset \mathbb R^n$ of $\bar x$, choose $r > 0$ and a neighbourhood $U_0 \subset \mathbb R^n$ of $\bar x$ such that
    \begin{align*}
        U_0 \subseteq B(\bar x, \delta^{-1} r), \qquad B(\bar x, r) \subseteq U_1 \cap U
    \end{align*}
    with $\delta = \frac{\hat\delta}{\Gamma p}$. 
    
    Given $x_0 \in U_0$, we prove by induction that the sequence $(x_k)_{k=0}^N$ generated by \eqref{eq:newton-iteration} for any $N \in \mathbb N$ (uniquely) exists and satisfies $x_k \in U_1 \cap U$ for all $k \in \{0, \ldots, N\}$: Since the claim is true for $N=0$, let it now also be true for $N \in \mathbb N$. In particular, $x_N \in U$ and there exists a (unique) $x_{N+1} \in \Phi(x_N) = \{ \bar x_N \}$ satisfying either
    \begin{align*}
        \| x_{N+1} - \bar x \| \leq \hat\delta \| x_0 - \bar x \| 
    \end{align*}
    if $N = 0$, or otherwise
    \begin{align*}
        \| x_{N+1} - \bar x \| \leq \hat\delta \Gamma p^N \| x_0 - \bar x \| 
    \end{align*}
    by virtue of the induction hypothesis. In either case, 
    \begin{align*}
        \| x_{N+1} - \bar x \| \leq \delta \| x_0 - \bar x \|
    \end{align*}
    as $p \in (0,1)$ and hence, $x_{N+1} \in U_1 \cap U$.

    By induction, the sequence $(x_k)_{k \in \mathbb N}$ exists (uniquely) and satisfies $x_k \in U_1 \cap U$ for all $k \in \mathbb N$; linear convergence follows immediately.
\end{proof}

\section{Conclusion}
\added[id=TC]{The Newton method for constrained nonconvex optimization, under common assumptions 
for the convex approximation,
admits a local vector Lyapunov function associated with constraint violation and suboptimality of the primal solutions. The construction of a strict Lyapunov function with quadratic bounds establishes local linear convergence and (thus exponential) stability. Our results provide a systems theoretical framework for nonconvex optimization algorithms and, in future work, globalization techniques.}

\appendix
\section{Strong Convexity implies Regularity}
Under the standing assumptions, let $(\bar x_k, \bar y_k) \in \Omega_k \times \mathbb R^m$ be the optimal primal-dual solution of \eqref{eq:quadratic-problem} and take $(x, y) \in \mathbb R^n \times \mathbb R^m$. By surjectivity of $\nabla g(x_k)$, there exists a solution $x' \in \Omega_k$ with
\begin{align*}
    \| x' - x \| \leq \chi \| g_k(x) \|
\end{align*}
which is \eqref{eq:metric-regularity-linear} with $z = 0$. Strong convexity of $\varphi_k$ yields
\begin{align*}
    \varphi_k(\bar x_k) \geq \varphi_k(x') + \langle \nabla L_k(x', y), \bar x_k - x' \rangle + \frac{\mu}{2} \| \bar x_k - x' \|^2
\end{align*}
and
\begin{align*}
    \varphi_k(x') \geq \varphi_k(\bar x_k) + \langle \nabla L_k(\bar x_k, \bar y_k), x' - \bar x_k \rangle + \frac{\mu}{2} \| x' - \bar x_k \|^2
\end{align*}
with
\begin{align*}
    \| \nabla L_k(x', y) - \nabla L_k(x, y) \| \leq \eta \| x' - x \|
\end{align*}
by Lipschitz continuity of $\nabla L_k(\cdot, y)$.
Combining the inequalities, we obtain
\begin{align*}
    \mu \| x' - \bar x_k \|^2 &\leq \langle \nabla L_k(x', y) - \nabla L_k(\bar x_k, \bar y_k), x' - \bar x_k \rangle
    \leq \| \nabla L_k(x', y) - \nabla \added[id=IK]{L_k}(\bar x_k, \bar y_k) \| \cdot \| x' - \bar x_k \|
\end{align*}
and hence, \added[id=IK]{$\| x' - \bar x_k \| \leq \mu^{-1} \| \nabla L_k(x', y) \|$ and}
\begin{multline*}
    \| x - \bar x_k \| \leq \| x' - \bar x_k \| + \| x' - x \| 
    \leq \mu^{-1} \| \nabla L_k(x', y) \| + \| x' - x \| \\
    \leq \mu^{-1} \| \nabla L_k(x, y) \| + \big( \frac{\eta}{\mu} + 1 \big) \| x' - x \| 
    \leq \mu^{-1} \| \nabla L_k(x, y) \| + \big( \frac{\eta}{\mu} + 1 \big) \chi \| g_k(x) \|
\end{multline*}
using the triangle inequality and $\nabla L_k(\bar x_k, \bar y_k) = 0$.

Strong convexity of $\varphi_k$ relative to any subspace $g_k^{-1}(b)$ with $b \in \mathbb R^m$ and surjectivity of $\nabla g(x_k)$ implies that the perturbed problem
\begin{align*}
    \min_{x \in \mathbb R^n} \varphi_k(x) - \langle b_0, x \rangle \quad \text{subject to $g_k(x) = b_1$}
\end{align*}
has an optimal primal-dual solution $(x_b, y_b) \in \mathbb R^n \times \mathbb R^m$ for any $(b_0, b_1) \in \mathbb R^n \times \mathbb R^m$ satisfying
\begin{align*}
    \| x_{b'} - x_b \| \leq \mu^{-1} \| b_0' - b_0 \| + \big( \frac{\eta}{\mu} + 1 \big) \added[id=IK]{\chi} \| b_1' - b_1 \|
\end{align*}
by virtue of the above subregularity property.

\bibliographystyle{IEEEtran}
\bibliography{references}

\end{document}